\documentclass[a4paper]{article}
\usepackage[ansinew]{inputenc}
\usepackage{enumerate}
\usepackage[reqno]{amsmath}
\usepackage{amsthm, amssymb}
\newtheorem{theorem}{Theorem}
\newtheorem{proposition}{Proposition}
\newtheorem{lemma}{Lemma}
\newtheorem{definition}{Definition}
\usepackage{graphicx}
\usepackage[left=2.5cm,right=2.5cm,top=1.5cm,%
bottom=3cm,includeheadfoot]{geometry}
\newcommand{\Leer}{\vspace{\baselineskip}}
\newcommand{\dis}{\displaystyle}
\newcommand{\N}{\mathbb{N}}

\newcommand{\C}{\mathbb{C}}
\newcommand{\R}{\mathbb{R}}
\newcommand{\Ls}{L}
\newcommand{\F}{\text{F}}

\newcommand{\sgn}{\operatorname{sgn}}

\allowdisplaybreaks
\begin{document}
\title{Partial sums of the M\"{o}bius function 
in arithmetic progressions assuming GRH}
\author{Karin Halupczok and Benjamin Suger}

\maketitle

\begin{abstract}
\noindent We consider Mertens' function in arithmetic progression,
\[
    M(x,q,a) := \underset{\substack{n\le x \\
    n\equiv a \text{ mod } q}}{\sum}\mu(n).
\]
Assuming the generalized Riemann hypothesis (GRH),
we show that the bound
\[
   M(x,q,a)\ll_{\varepsilon}
    \sqrt{x}\exp{\left((\log x)^{3/5}(\log\log x)^{16/5 +
      \varepsilon}\right)}
\]
holds uniform for all $q\le \exp\Big(\frac{\log 2}{2}
   \Big\lfloor (\log x)^{3/5}(\log\log x)^{11/5}\Big\rfloor\Big)$,
$\gcd(a,q)=1$ and all $\varepsilon >0$. The implicit constant 
is depending only on $\varepsilon$.
For the proof, a former method of K.\ Soundararajan 
is extended to $L$-series.
\end{abstract}



\tableofcontents


\newpage
\section{Introduction}
\label{sec1}

Mertens' function is defined by
\[M(x) := \underset{n\le x}{\sum}\mu(n). \]
It is well known that
$M(x)=\mathsf{o}(x^{1/2 + \varepsilon})$ is equivalent 
to Riemann's hypothesis.

When assuming Riemann's hypothesis for $\zeta$,
one can give even sharper bounds for $M(x)$, see \cite{Lan},
   \cite{Tit}, \cite{MaMo}, \cite{Sou}, \cite{RoBa}:

In \cite{MaMo}, H.\ Maier and H.\ L.\ Montgomery proved the bound
\[
   M(x)\ll x^{1/2}\exp\Big(c(\log x)^{39/61}\Big) \text{ for a }c>0.
\]
In \cite{Sou}, K.\ Soundararajan improved the bound by showing
\[
    M(x)\ll x^{1/2}\exp\Big((\log x)^{1/2}(\log\log x)^{14}\Big).
\]
In \cite{RoBa}, A.\ de Roton und M.\ Balazard refine
the result of K.\ Soundararajan and show
\[
    M(x) \ll_{\varepsilon} \sqrt{x}\exp\Big((\log
    x)^{1/2}(\log\log x)^{5/2 + \varepsilon}\Big),
\]  
which is the best bound up to date.

In this paper we generalize the method of K.\ Soundararajan to
provide a bound for Mertens' function in arithmetic progression,
\[
   M(x,q,a) := \underset{\substack{n\le x \\
    n\equiv a \text{ mod } q}}{\sum}\mu(n).
\]
Note that the trivial bound is $\leq x/q$, so bounds smaller than
$x^{1/2+\varepsilon}$ are nontrivial if $q\leq x^{1/2-\varepsilon}$.

We adapt the method of K.\ Soundararajan 
resp.\ the modification of A.\ de Roton and M.\ Balazard 
in such a way, that it remains applicable for
Dirichlet $L$-series. We obtain the following nontrivial 
upper bound assuming Riemann's hypothesis for all Dirichlet 
$L$-series $L(s,\chi)$ with $\chi$ mod $q$ and all moduli $q$
in question (GRH for short):

\begin{theorem}
  \label{t}
Assuming GRH, the bound
\[M(x,q,a)\ll_{\varepsilon}
\sqrt{x}\exp{\left((\log x)^{3/5}(\log\log x)^{16/5 +
      \varepsilon}\right)} \] 
holds uniform for all 
$q\le \exp\Big(\frac{\log 2}{2}
   \Big\lfloor (\log x)^{3/5}(\log\log x)^{11/5}\Big\rfloor\Big)$,
$\gcd(a,q)=1$ and all $\varepsilon >0$ with an implicit constant
depending only on $\varepsilon$.
\end{theorem}

With this theorem, we extend the results of \cite{Sou} resp.\
\cite{RoBa} to a Siegel-Walfisz-type result. The obtained bound
is weaker than the one of \cite{Sou} resp.\ \cite{RoBa},
but still sharper than the one of \cite{MaMo}.

The method we use is as follows.
We expand the M\"obius sum $M(x,q,a)$ using Dirichlet characters,
\begin{align*}
    M(x,q,a)&= \frac{1}{\varphi(q)} \sum_{\chi(q)} \bar{\chi}(a)
    \sum_{n\leq x}\chi(n)\mu(n) \\
     &= \frac{1}{\varphi(q)} \sum_{\chi(q)} \bar{\chi}(a) A(x,\chi,q)
     + {\rm O}(\log x),
\end{align*}
using Perron's formula with integrals
\[
   A(x,q,\chi)= \frac{1}{2\pi i} \int_{1+1/(\log x)
     -i2^{K}}^{1+1/(\log x) + i2^{K}} \frac{x^{s}}{L(s,\chi)s} ds,
   \quad K:= \Big\lfloor \frac{\log x}{\log 2} \Big\rfloor.
\]
With this, bounds for $L(s,\chi)$ are needed. Considering the
principal character $\chi_{0}$ mod $q$, the formula
\[
   L(s,\chi_{0}) = \zeta(s) \prod_{p|q} \Big(1-\frac{1}{p^{s}}\Big)
\]
shows that already the sharper bound of \cite{Sou}/\cite{RoBa}
applies (see the proof of Lemma~\ref{l3}).
So the main work is to consider nonprincipal characters.

Like in \cite{Sou}/\cite{RoBa}, the main steps are then some
propositions aiming to bound $L(s,\chi)$ to obtain
an upper estimate for $A(x,q,\chi)$. They are given in Sections~\ref{sec7}
and \ref{sec8} and are resulting from the propositions in the former Sections
~\ref{sec2} and \ref{sec4},\ref{sec5},\ref{sec6}.

Most of these propositions are stated for primitive 
characters. If necessary, 
results for nonprimitive characters $\chi\neq \chi_{0}$ 
are derived by reduction to a primitive character that induces 
$\chi$.

The main idea in \cite{Sou}, namely the concept of $V$-typical
ordinates, is extended to a version which allows one to work
also with $L$-series. We give the adapted definition in Section~\ref{sec3}.

As one important step, we show in Section \ref{sec4} that there
\textit{are} actually $V$-typical ordinates, see Proposition~\ref{Prop9}.

In Section~\ref{sec5}, it is shown that short 
intervals containing an unusual number of ordinates
of nontrivial $L$-zeros mod $q$ do not appear too often,
even uniformly for all $q$ up to the given bound (Proposition \ref{Prop8}),
so the $V$-untypical ordinates are small in number (Proposition \ref{Prop10}).
In the case of $\zeta$, this has been the breakthrough in
Soundararajan's paper \cite{Sou}. 

The resulting bound and the range for $q$ in Theorem~\ref{t}
is then obtained by optimizing the bounds 
for $A(x,q,\chi)$ in Section~\ref{sec9}.
The elementary Proposition \ref{Prop20} plays an intrinsic r\^ole for this.


\vspace{\baselineskip}
A remark on notations used in this paper:

We mark all Propositions that assume the generalized Riemann
hypothesis by the symbol (GRH). 
We stress that 
all implicit constants are absolute unless otherwise indicated.


\section{List of tools}
\label{sec2}

In this section, we give a collection of the tools used in the proof.

The first proposition gives an approximation of the characteristic
function of a given interval:
\newcommand{\eins}{\mathbf{1}}

\begin{proposition}
\label{Prop1}
Let $h>0$, $\Delta \ge 1$.

Let $\eins_{\lbrack -h,h \rbrack}$ be the characteristic function of
the interval $\lbrack -h, h \rbrack$.

There are even, entire functions $F_{+}$ and $F_{-}$
depending on $h$ and $\Delta$, being real on the real axis
and such that the following properties hold:

\begin{enumerate}
\item $ \forall u\in \R:~ F_{-}(u)\le \eins_{\lbrack -h,h \rbrack}(u) 
        \le F_{+}(u)$,
\item $\underset{-\infty}{\overset{\infty}{\int}}|F_{\pm}(u) -
  \eins_{\lbrack -h,h \rbrack}(u)|du = 1/\Delta$  and
  $\hat{F}_{\pm}(0) = 2h \pm 1/\Delta$,
\item $\hat{F}_{\pm}$ is realvalued and even, and we have
	 $\hat{F}_{\pm}(x) = 0$ for all $|x|\ge\Delta$ and 
      $|x\hat{F}_{\pm}(x)|\le 2$ for all $x\in\R$,
\item for $z\in\C$ with $|z|\ge \max\{2h,1\}$ we have
	\[|F_{\pm}(z)|\ll \frac{\exp(2\pi |\Im z|\Delta)}{(\Delta
          |z|)^{2}}.\] 
\end{enumerate}
\end{proposition}

The proof uses Beurling's Approximation of the signum
function
\[
  \operatorname{sgn}(x):=
  \begin{cases}
    x/|x|, & x\neq 0,\\
    0,     & x=0.
  \end{cases}
\]
Let 
$K(z):=\Big(\frac{\sin(\pi z)}{\pi z}\Big)^{2}$
and
$H(z)=K(z)\Big(\sum_{n=-\infty}^{\infty}
\frac{\operatorname{sgn}(n)}{(z-n)^{2}} + \frac{2}{z} \Big),$
then it can be shown that the functions
\[
  F_{\pm}(z):= \frac{1}{2} (H(\Delta(z+h)) \pm K(\Delta(z+h)) 
     + H(\Delta (h-z)) \pm K(\Delta (h-z)) )
\]
have the properties asserted in Proposition \ref{Prop1}.
This can be seen as in \cite{Sel} and \cite{Val}, see also \cite{RoBa},
we just give the proof of part 4.\ in more detail:

For this, let $z=x+iy$ with $x,y\in\R$ and $|z|\ge\max\{2h,1\}$. Since 
$\F_{\pm}$ are even, consider only nonnegative $x$.
Using $\sin(z) \ll e^{|\Im(z)|}$
and $\Im(\Delta(z+h)) = -\Im(\Delta(h-z)) = \Delta\Im(z)$, we get the 
desired bound for $K(\Delta(z+h)) \pm K(\Delta(h-z))$ 
since $|z \pm h| = |z|\left|1 \pm \frac{h}{z}\right| \ge
|z|\left(1-\frac{h}{|z|}\right) \ge \frac12 |z|$.

To estimate $H(\Delta(z+h)) + H(\Delta(h-z))$
we use the identities
\begin{align}
\left(\frac{\pi}{\sin(\pi z)}\right)^2 =
\sum_{n=-\infty}^{\infty}\frac{1}{(z-n)^2}, 
&\text{ converging on every compact subset of }\C\setminus\mathbb{Z},
\label{sin}\\
\underset{n=0}{\overset{\infty}{\sum}}\frac{1}{(z+n)(z+n+1)} =
\frac{1}{z}, &\text{ converging absolutely for } z\in\C\setminus-\N_0.
\label{einsdurch}
\end{align}

Consider $H(\Delta(z+h))$ and $H(\Delta(h-z))$ separately. By
\eqref{sin}, we have
\begin{align*}
H(\Delta(z+h)) &= \left(\frac{\sin(\pi\Delta(z+h))}{\pi}\right)^2\left(\sum_{n=-\infty}^{\infty}\frac{\sgn(n)}{(\Delta(z+h)-n)^2} +\frac{2}{\Delta(z+h)} \right)\\
        &=1 +
        \left(\frac{\sin(\pi\Delta(z+h))}{\pi}\right)^2\left(-2\sum_{n=1}^{\infty}\frac{1}{(\Delta(z+h)+n)^2}
          - \frac{1}{(\Delta(z+h))^2}  +\frac{2}{\Delta(z+h)} \right),\\ 
\end{align*}
and \eqref{einsdurch} gives for the negative of the last term in large
brackets the expression
\begin{align*}
&\sum_{n=0}^{\infty}\left(\frac{1}{(\Delta(z+h)+n)^2} + \frac{1}{(\Delta(z+h)+n+1)^2}\right) - \sum_{n=0}^{\infty}\frac{2}{(\Delta(z+h)+n)(\Delta(z+h)+n+1)} \\
&= \sum_{n=0}^{\infty}\left(\frac{1}{(\Delta(z+h)+n)}
  - \frac{1}{(\Delta(z+h)+n+1)}\right)^2 
=\sum_{n=0}^{\infty}\frac{1}{(\Delta(z+h)+n)^{2}(\Delta(z+h)+n+1)^2} \\
&\le \frac{1}{(\Delta(x + h + |y| ))^2}\sum_{n=0}^{\infty}\frac{1}{(\Delta(x + h + |y| )+n )(\Delta(x + h + |y| )+n+1)} \\
&=\frac{1}{(\Delta(x + h + |y| ))^3} \ll \frac{1}{|\Delta(z + h )|^3}
\ll \frac{1}{|\Delta z|^{3}}.
\end{align*}
Analogously, we get
\begin{align*}
H(\Delta(h-z)) &= \left(\frac{\sin(\pi\Delta(h-z))}{\pi}\right)^2\left(\underset{-\infty}{\overset{\infty}{\sum}}\frac{\sgn(n)}{(\Delta(h-z)-n)^2} +\frac{2}{\Delta(h-z)} \right)\\ 
&= -1 + \left(\frac{\sin(\pi\Delta(z-h))}{\pi}\right)^2\left(\frac{1}{(\Delta(z-h))^2} +2\underset{1}{\overset{\infty}{\sum}}\frac{1}{(\Delta(z-h)+n)^2}  -\frac{2}{\Delta(z-h)} \right).
\end{align*}
If $\Re(z)>h$, the treatment of the last term in large brackets is as before.

So let $\Re(z)\le h$. Due to $|z|\ge 2h$, we have $|y|=|\Im(z)|>h$, so
$z\not\in\R$ and $|\Im(z)|\geq |\Re(z)|$. Again (\ref{einsdurch})
gives for the last term in large brackets the expression
\begin{align*}
& \sum_{n=0}^{\infty} \frac{1}{(\Delta(z-h)+n)^{2}(\Delta(z-h)+n+1)^{2}}\\
&\ll \sum_{0\leq n \leq \Delta h}
 \frac{1}{|(\Delta(x-h)+n| + \Delta|y|)^{2}(|\Delta(x-h)+n+1| +
     \Delta|y|)^{2}} \\
&+ \sum_{n>\Delta h}
 \frac{1}{|(\Delta(x-h)+n| + \Delta|y|)^{2}(|\Delta(x-h)+n+1| + \Delta|y|)^{2}}\\
&\ll\frac{\max\{\Delta h,1\}}{|\Delta y|^4} + \sum_{n=0}^{\infty}
\frac{1}{(\Delta|y|+n)^{2}(\Delta|y|+n+1)^{2}} 
\ll \frac{1}{|\Delta y|^3}  \ll \frac{1}{|\Delta z|^3}.
\end{align*}
Summing up we obtain
\[H(\Delta(z+h)) + H(\Delta(h-z)) \ll
\frac{e^{2\pi\Delta|\Im(z)|}}{(\Delta|z|)^3} \]
and the desired bound for $|z|\ge\max\{2h,1\}$. \qed

\Leer
We will make use of the following explicit formula
for the functions $F_{\pm}$.

\begin{proposition}
\label{Prop2}
(GRH) Let $\chi$ be a primitive character mod $q$.
Let $t>0$, $\Delta \ge 1$, $h>0$, and $F_{\pm}$ 
the functions from Proposition \ref{Prop1}.
Then we have
\begin{multline*}
    \underset{\rho = \frac{1}{2} + i\gamma}{\sum}F_{\pm}(\gamma -t) =
    \frac{1}{2\pi}\hat{F}_{\pm}(0) \log \frac{q}{\pi}  +
    \frac{1}{2\pi}\underset{-\infty}{\overset{\infty}{\int}}
       F_{\pm}(u-t)\Re\frac{\Gamma'}{\Gamma}\left(\frac{\frac{1}{2} 
        + iu + \mathsf{a}}{2}\right)du \\
    - \frac{1}{\pi}\Re\underset{n\in\N}{\sum}
        \frac{\Lambda(n)\chi(n)}{n^{\frac12
        +it}}\hat{F}_{\pm}\left(\frac{\log n}{2\pi}\right).
\end{multline*}
Here the sum on the left hand side runs through all zeros of
$\Ls(s,\chi)$ in the strip $0\leq \sigma\leq 1$ with relevant 
multiplicity, and where we have set
\begin{equation}
\label{adef}
\mathsf{a} := \mathsf{a}(\chi) := 
\begin{cases} 0, & \text{ if } \chi(-1) = 1, \\
                     1, & \text{ if } \chi(-1) = -1.
\end{cases}
\end{equation}
\end{proposition}

The proof can be established in the same way as Theorem 5.12, p.\ 108,
in the book \cite{IwKo} of Iwaniec and Kowalski. It uses the Mellin
transform, the explicit formula for $\frac{L'}{L}(s,\chi)$ and the
residue theorem, where one has to take care of the trivial zero of 
$L(s,\chi)$ at $s=0$ if $\chi(-1)=1$.

\Leer
An estimate of the integral in Proposition \ref{Prop2} gives 
the next proposition:

\begin{proposition}
\label{Prop3}
Let $t\ge 25$, $\Delta\ge 1$, $0<h\le \sqrt{t}$, $F_{\pm}$ as in
Proposition \ref{Prop1}, $\chi$ a character mod $q$.
Then it holds that
\[ \underset{-\infty}{\overset{\infty}{\int}}F_{\pm}(u-t)
\Re\frac{\Gamma'}{\Gamma}\left(\frac{1}{4} + \frac{\mathsf{a} +
    it}{2}\right) du = \left(2h \pm
  \frac{1}{\Delta}\right)\log\frac{t}{2} + {\rm O}(1),  \]
where $\mathsf{a}$ is defined in (\ref{adef}).
\end{proposition}

The proof can be obtained as in \cite{GoGo}.
It uses Stirling's formula and the properties of $F_{\pm}$ from
Proposition \ref{Prop1} after splitting the integral at
$t-4\sqrt{t}$ and $t+4\sqrt{t}$.

\Leer
We make also use of the following result of Maier and Montgomery 
in \cite{MaMo} concerning moments of Dirichlet polynomials:

\begin{proposition}
\label{Prop4}
Consider a Dirichlet polynomial
$P(s) = \sum_{p\le N} a(p)p^{-s}$. For $T\geq 3$
and $\alpha \in\R$ let
$s_1 , ... , s_R \in\C$ with $1\le |\Im (s_i - s_j)|\le T$ for
$i\not= j$, and $\Re s_i \ge \alpha$ for $ 1\le i \le R$.

Then, for every positive integer $k$ with $N^{k}\le T$, it holds that
\[\underset{r=1}{\overset{R}{\sum}}|P(s_r)|^{2k} \ll T (\log T)^{2} k!
\Big(\sum_{p\le N}|a(p)|^{2}p^{-2\alpha}\Big)^{k}.\]
\end{proposition}



Our result relies further on the estimate in the 
following proposition.

\begin{proposition}
\label{Prop5}
Let $T \ge e^{e^{33}}$,
$\left(\log \log T\right)^{2}\le V \le \frac{\log T}{\log\log T}$, 
$\eta = \frac{1}{\log V}$ and
$k = \left\lfloor\frac{2V}{3(1 + \eta)}\right\rfloor$.

Then we have
\[ k\left(\log(k\log\log T) - 2\log(\eta V)  \right) \le -\frac23
V\log\frac{V}{\log\log T} + \frac43 V\log\log V + \frac23 V. \]
\end{proposition}

The proof is completely analogous to the elementary proof in 
\cite{RoBa}, there Proposition 14 on page 11 and 12.

\Leer
Now using Proposition \ref{Prop2},
we can give an upper and lower bound for the number of zeros
in a certain region around ordinate $t$.

\begin{proposition}
\label{Prop6}
(GRH) Let $t\ge 25$, $\Delta \ge 2$, $0< h \le\sqrt{t}$ and
$\chi$ be a primitive character mod $q$. Then
\begin{equation*}
 -\frac{\log(qt)}{2\pi\Delta} - \frac{1}{\pi}\Re \underset{p\le
   e^{2\pi\Delta}}{\sum}\frac{\chi(p) \log(p)}{p^{\frac{1}{2}
     +it}}\hat{F}_{-}\left(\frac{\log p}{2\pi} \right) +
 {\rm O}(\log \Delta) \le  
 N(t+ h,\chi) - N( t - h , \chi) - \frac{h}{\pi}\log\frac{qt}{2\pi}
 \end{equation*}
and
\begin{equation*}
N(t+ h,\chi) - N( t - h , \chi) -
\frac{h}{\pi}\log\frac{qt}{2\pi}\le\frac{\log(qt)}{2\pi\Delta} -
\frac{1}{\pi}\Re \underset{p\le e^{2\pi\Delta}}{\sum}\frac{\chi(p)
  \log(p)}{p^{\frac{1}{2} +it}}\hat{F}_{+}\left(\frac{\log p}{2\pi}
\right) + {\rm O}(\log \Delta).
\end{equation*}
\end{proposition}

\textbf{Proof:}

We only show the upper bound, the lower bound estimate can be done
in a complete analogous way.

We use the functions of Proposition \ref{Prop1} and the results
from Propositions \ref{Prop2} and \ref{Prop3}, we see analogously to
\cite{RoBa} (there Proposition 15 from page 12 on):
\begin{equation*}
N(t+h,\chi) - N(t-h,\chi) \le \left(2h +
  \frac{1}{\Delta}\right)\frac{1}{2\pi}\log\frac{qt}{2\pi} +
{\rm O}(1) - \frac{1}{\pi}\Re\underset{n\le
  e^{2\pi\Delta}}{\sum}\frac{\Lambda(n)\chi(n)}{n^{\frac{1}{2} +
    it}}\hat{F}_{+}\left(\frac{\log n}{2\pi}\right).
\end{equation*}

Here
\begin{align*}
\frac{1}{\pi}\Re\underset{n\le
  e^{2\pi\Delta}}{\sum}\frac{\Lambda(n)\chi(n)}{n^{\frac{1}{2} +
    it}}\hat{F}_{+}\left(\frac{\log n}{2\pi}\right) &=  
\frac{1}{\pi}\Re\underset{p\le e^{2\pi\Delta}}{\sum}\frac{\log p ~
  \chi(p)}{p^{\frac{1}{2} + it}}\hat{F}_{+}\left(\frac{\log
    p}{2\pi}\right) \\ &+\frac{1}{\pi}\Re\underset{p\le
  e^{\pi\Delta}}{\sum}\frac{\log p ~ \chi(p)^2}{p^{1 +
    2it}}\hat{F}_{+}\left(\frac{\log p}{\pi}\right) +
{\rm O}(1)\\  
&=\frac{1}{\pi}\Re\underset{p\le e^{2\pi\Delta}}{\sum}\frac{\log p ~
  \chi(p)}{p^{\frac{1}{2} + it}}\hat{F}_{+}\left(\frac{\log
    p}{2\pi}\right) + {\rm O}(\log\Delta),
\end{align*}
and this finishes the proof. \qed


\section{V-typical ordinates}
\label{sec3}

The method of Soundararajan in \cite{Sou}
relies on the notion of $V$-typical ordinates.
We modify this definition for our purposes and define
$V_{(\delta, \chi , \text{q})}$-typical ordinates as follows.

\begin{definition}
\label{Def1} ($V_{(\delta, \chi , q)}$-typical).

Let $q\in\N$ and $\chi$ a character mod $q$.
If $\chi$ is nonprincipal, let it be induced by
$\chi_{1}$ mod $q_1$, 
let $T>e$ and $0< \delta \le 1$. 

Let $V\in \Bigl\lbrack (\log\log T)^2, \dis
\frac{\log T}{\log\log T} \Bigr\rbrack$.

An ordinate $t\in \lbrack \text{T, 2T} \rbrack$ is called 
\underline{$V_{(\delta, \chi , q)}$-typical} of order $T$, if the following 
properties hold:

\begin{enumerate}[(i)]
\item $\forall \sigma\ge\frac{1}{2}$: $\dis\Bigl|
\sum_{n \le x} 
   \frac{\chi_{1}(n) \Lambda(n)}{n^{\sigma + it}\log n} 
   \frac{\log\bigl(\frac{x}{n}\bigr)}{\log x} \Bigr|
     \le 2V$ with $x=T^{\frac{1}{V}}$, 
\item $\forall t' \in ( t-1 , t+1 )$:\\ $N(t' + h, \chi) - N(t' - h, \chi) \le
 (1+\delta) V$ with  $
 h=\dis\frac{\delta \pi V}{\log(q_1T)}$ and $\lbrack t'-h,t'+h \rbrack
 \subseteq \lbrack t-1,t+1 \rbrack $,
\item
 $ \forall t'\in (t-1,t+1)$:\\ $N(t' + h, \chi ) - N(t'-h, \chi) \le V $
   with $ h=\dis\frac{\pi V}{\log V \log (q_1T)}$ and 
   $\lbrack t'-h,t'+h \rbrack \subseteq \lbrack t-1,t+1 \rbrack $.
\end{enumerate}

If at least one of the three properties does not hold, we call $t$ a
\underline{$V_{(\delta, \chi , q)}$-untypical} ordinate of order $T$.
\end{definition}

In what follows, the meaning of $\chi$,
$q$ and $\delta$ is often clear from the context, then we will
write simply $V$-typical instead of 
$V_{(\delta, \chi , q)}$-typical of order $T$.



\section{$V$ such that
  all $t\in\lbrack T,2T\rbrack$ are $V$-typical}
\label{sec4}

\Leer
\begin{proposition}
\label{Prop7}
Let $t$ be sufficiently large and let $0<h \le \sqrt{t}$, let $\chi$ 
be a primitive character mod $q$. Then
\[
  \left| N(t+h,\chi) - N(t-h,\chi) -
  \frac{h}{\pi}\log\frac{qt}{2\pi} \right| \le
\frac{\log(qt)}{2\log\log (qt)} + \left( \frac{1}{2} +
  {\rm o}\left(1\right) \right)\frac{\log (qt) \log\log\log
  (qt)}{(\log\log (qt))^{2}} \text{ for } t\to \infty.
\]
\end{proposition}

\textbf{Proof:}
As in \cite{RoBa}, 
we estimate the sum of Proposition \ref{Prop6}
as follows:
\begin{equation}
\label{oben}
\left|\frac{1}{\pi}\Re\underset{p\le e^{2\pi\Delta}}{\sum}\frac{\log p
    \,\chi(p)}{p^{\frac{1}{2} + it}}\hat{F}_{+}\left(\frac{\log
      p}{2\pi}\right)\right| \ll \sum_{p\le  e^{2\pi\Delta}}
    \frac{1}{\sqrt{p}} \ll \frac{e^{\pi\Delta}}{\Delta}.
\end{equation}

Now set $\Delta = \dis\frac{1}{\pi}\log{\frac{\log (qt)}{\log\log (qt)}}$.
By estimate (\ref{oben}), we obtain
\begin{align*}
\Bigl| N(t+h,\chi) - &N(t-h,\chi) - \frac{h}{\pi}\log\frac{qt}{2\pi} \Bigr| \\ 
&\le\frac{\log(qt)}{2(\log\log(qt) - \log\log\log(qt))} +
{\rm O}\Big(\frac{\frac{\log(qt)}{\log\log(qt)}}{\log\log(qt) -
    \log\log\log(qt)}  \Big) \\ 
&=\frac{\log(qt)}{2\log\log(qt)}\underset{k=0}{\overset{\infty}{\sum}}\left(
  \frac{\log\log\log(qt)}{\log\log(qt)}
\right)^{k}+{\rm O}\left(\frac{\log(qt)}{(\log\log(qt))^{2}}\right)\\ 
&=\frac{\log(qt)}{2\log\log(qt)} + \frac{\log(qt)
  \log\log\log(qt)}{2(\log\log(qt))^2}(1+{\rm o}(1))
\end{align*}
with an ${\rm o}(1)$-term not depending on $q$, more precise, it is
${\rm O}((\log\log\log t)^{-1} )$. \qed

\Leer
\begin{proposition}\label{Prop9}
Let $\chi$ be a character mod $q$, $q_1$ 
be the conductor of $\chi$ and $0<\delta\le1$.
Further let $T$ be sufficiently large, at least
$T\ge \max\{q^{2},e^{e^9}\}$, and let $V$ be such that
\[
  \frac{3}{4} + \frac{\log\log\log T}{\log\log T}\le V \frac{\log\log
  T}{\log T} \le 1
\]
holds.
Then all ordinates $t\in\lbrack T,2T \rbrack$ are
$V$-typical of order $T$.
\end{proposition}

As a consequence of this proposition, we conclude that 
$V$-typical ordinates exist.

\textbf{Proof:}

We have to verify properties (i), (ii) and (iii) from Definition~
\ref{Def1}.

\Leer
Ad (i):

Let $f(u) := \underset{2\le n\le u }{\sum}
  \frac{\Lambda(n) \chi_{1}(n)}{\sqrt{n} \log n}$, $u\ge 2$.
Then (see \cite{RoBa}, page 16):
\begin{equation*}
|f(u)| \le \underset{2\le n\le u}{\sum}
  \frac{\Lambda(n)}{\sqrt{n}\log n} \ll \frac{\sqrt{u}}{\log u},
\end{equation*}
and from this we obtain
\[ \underset{n\le x}{\sum}\frac{\Lambda (n) \chi_{1}(n)}{\sqrt{n} \log
  n}\log\frac{x}{n} = \int_{1}^{x}\frac{f(u)}{u}du
\ll \frac{\sqrt{x}}{\log x}. \]
Since $ x=T^{\frac{1}{V}} \le T^{\frac{4\log\log T}{3\log T}} \le (\log T)^2$,
we have
\begin{equation*}
  \left|\underset{n \le x}{\sum}
  \frac{\chi_{1}(n) \Lambda(n)}{n^{\sigma + it}\log n} 
  \frac{\log\left(\frac{x}{n}\right)}{\log x} \right| 
  \ll \frac{\sqrt{x}}{(\log x)^2} 
  \ll \frac{\log T}{(\log\log T )^2} = {\rm o} (V).
\end{equation*}

\Leer
Ad (ii):

Let $t'\in \lbrack t-1,t+1\rbrack$ and $h= \frac{\delta \pi
  V}{\log(q_1T)} $.
Since $h=\frac{\delta\pi V}{\log (q_1T)}\le 
\pi V \le \log T < \sqrt{T}$, we can apply Proposition
\ref{Prop7} on the primitive character $\chi_1$ mod $q_1$
that induces $\chi$. We obtain, using $q^{2}\le T$, that
\begin{align*}
N(t'+h,\chi) &- N(t'-h, \chi) \\
&\le \frac{h}{\pi}\log\frac{q_1t'}{2\pi} +
\frac{\log(q_1t')}{2\log\log (q_1t')} + \Big( \frac{1}{2} +
  {\rm o}(1) \Big)\frac{\log (q_1t') \log\log\log
  (q_1 t')}{(\log\log (q_1 t'))^{2}} \\ 
&\le\frac{h}{\pi}\log\frac{q_1 T}{\pi} + \frac{\log(2qT)}{2\log\log T}
+ \Big( \frac{1}{2} + {\rm o}(1) \Big)\frac{\log (2qT)
  \log\log\log T}{(\log\log T)^{2}} \\ 
&\le \delta V + \frac{\log T^{3/2}}{2\log\log T} +
+ \Big(\frac{1}{2} + {\rm o}(1) \Big)\frac{\log T^{3/2}
  \log\log\log T}{(\log\log T)^{2}} \\ 
&=\delta V + \frac{3\log T}{4\log\log T} +  \Big( \frac{3}{4} +
  {\rm o}(1) \Big)\frac{\log T \log\log\log T}{(\log\log T)^{2}}
\\ 
&\le \delta V + \frac{3\log T}{4\log\log T} + \frac{\log T
  \log\log\log T}{(\log\log T)^{2}} \\ 
&\dis\le(1 + \delta)V .
\end{align*}

\Leer
Ad (iii):

Let $t'\in \lbrack t-1,t+1\rbrack$ and 
$h= \frac{\pi V}{\log V \log(q_1T)} $, then
\begin{align*}
N(t+h,\chi) &- N(t-h, \chi) \\ 
&\le \frac{h}{\pi}\log\frac{q_1t'}{2\pi} +
\frac{\log(q_1t')}{2\log\log (q_1t')} + \Big( \frac{1}{2} +
  {\rm o}(1) \Big)\frac{\log (q_1t') \log\log\log
  (q_1t')}{(\log\log (q_1t'))^{2}} \text{ by Prop. }\ref{Prop7} \\ 
&\le \frac{V}{\log V} + \frac{3\log
  T}{4\log\log T} + \Big( \frac{3}{4} + {\rm o}(1)
\Big)\frac{\log T \log\log\log T}{(\log\log T)^{2}} 
\text{ analoguously to (ii)}\\ 
&= \frac{3\log T}{4\log\log T} + \Big( \frac{3}{4} + {\rm o}(1)
\Big)\frac{\log T \log\log\log T}{(\log\log T)^{2}} \\ 
&\le \frac{3\log T}{4\log\log T} + \frac{\log T \log\log\log
  T}{(\log\log T)^{2}} \le V. 
\end{align*}
\qed



\section{The number of $V$-untypical, well separated
ordinates}
\label{sec5}

\Leer
\begin{proposition}\label{Prop8}
Let $\chi\not=\chi_0$ be a character mod $q$ and
$q_1$ be the conductor of $\chi$.
Further let
\begin{enumerate}
\item $T$ be large, at least $T\ge q^2$,
\item $0<h\le \sqrt{T}$,
\item $(\log\log T)^{2} \le V \le \dis\frac{\log T}{\log\log T}$,
\item $T\le t_1 < t_2 < \dots < t_R \le 2T$ and
              $t_{r+1} - t_r \ge 1$ for $1\le r < R$,
\item $N(t_r +h,\chi)- N(t_r -h,\chi) - 
          \dis\frac{h}{\pi}\log\frac{q_1t_r}{2\pi} \ge
           V + {\rm O}(1)$ 
          for $1\le r \le R$.
\end{enumerate}	
Then
\[ 
R \ll T \exp\Big(-\frac23 V\log\frac{V}{\log\log T}
   + \frac43 V\log\log V + {\rm O}(V)  \Big).
\]
\end{proposition}

\textbf{Proof:}

If $q_1=q$, then $\chi$ is primitive. If $q_1<q$, then $\chi$ 
is induced by a primitive character $\chi_1$ mod $q_1$, and we have
\[N(t,\chi) = N(t,\chi_1). \]

Therefore we can apply the results from Proposition \ref{Prop6} for
$\chi_1$ and $q_1$.
By the estimate from Proposition \ref{Prop6} we obtain
\begin{align*}
V + {\rm O}(1)&\le N(t_r +h,\chi_1)- N(t_r -h,\chi_1) -
\frac{h}{\pi}\log\frac{q_1t_r}{2\pi} \\ 
&\le \frac{\log(2qT)}{2\pi\Delta} + \left|\frac{1}{\pi}\underset{p\le
    e^{2\pi\Delta}}{\sum}\frac{\chi(p)\log p}{p^{\frac{1}{2} +
      it_r}}\hat{\F}_{+}\left(\frac{\log p}{2\pi}\right)\right| +
{\rm O}(\log\Delta), \quad\Delta\geq 2.
\end{align*}

If we define $a(p):=\frac{\chi(p)\log
  p}{\pi}\hat{F}_{+}\left(\frac{\log p}{2\pi}\right)$, we have: 
\[\left|\underset{p\le
    e^{2\pi\Delta}}{\sum}\frac{a(p)}{p^{\frac{1}{2}+ it_r}}\right| \ge
V - \frac{\log(2qT)}{2\pi\Delta} + {\rm O}(\log\Delta) + {\rm O}(1),\]
where $|a(p)|\le 4$ holds by Proposition \ref{Prop1}.

Let
\[\eta = \frac{1}{\log V} \text{ and } 
 \Delta = \frac{(1 + \eta) \log (qT)}{2\pi V}. \]

Then we have
\[ \exp(2\pi\Delta) 
   = (qT)^{\frac{1+\eta}{V}} \le 
  T^{\frac{3(1+\eta)}{2V}}  \text{ since } q\le\sqrt{T},
\]
hence
\[\log\Delta \ll \log\log T \le \sqrt{V}. \]

We obtain
\begin{align*}
 V - \frac{\log(2qT)}{2\pi\Delta} + {\rm O}(\log\Delta) 
  + {\rm O}(1)&= V - \frac{V\log(2qT)}{\left(1 + \eta\right) \log(qT)} +
 {\rm O}\left(\sqrt{V}\right) \\ 
 &\ge \frac{\eta V}{1+\eta}
 - \frac{\log2}{(1 + \eta)\log\log T} +
 {\rm O}\left(\sqrt{V}\right) \ge \frac{1}{2}\eta V.
\end{align*}

So we have
\[ \Bigg| \underset{p\le
    e^{2\pi\Delta}}{\sum}\frac{a(p)}{p^{\frac{1}{2}+ it_r}}  \Bigg|
\ge \frac{1}{2}\eta V  
   \text{ for } 1\le r \le R. \]

Let $k\in\N$ with
$k\le\left\lfloor\frac{2V}{3(1+\eta)}\right\rfloor$. Then we can apply
Proposition \ref{Prop4} with $N=(qT)^{(1+\eta)/V}$
since $(qT)^{k\frac{1+\eta}{V}} \le T^{k\frac{3(1+\eta)}{2V}} \le T$
for $q^{2}\le T$.
 
Raising to the $2k$-th power and summing over all $r=1,\dots ,R$,
applying Proposition \ref{Prop4} for $\alpha=\frac12$ and
$N=\left\lfloor(qT)^{\frac{1+\eta}{V}}\right\rfloor$,
we obtain analogously to \cite{RoBa} (page 15):
\begin{equation*}
R\Big( \frac{\eta V}{2} \Big)^{2k} \le
\underset{r=1}{\overset{R}{\sum}}\Bigg|\underset{p\le (qT)^{\frac{1 +
        \eta}{V}}}{\sum} \frac{a(p)}{p^{\frac{1}{2} + it_r}}
\Bigg|^{2k} \ll T (\log T)^2 (C k \log\log T)^k 
\end{equation*}
with an absolute constant $C>0$. So we have by now
\[R \ll T (\log T)^2 (4C)^k \Big( \frac{k \log\log T}{\eta^2 V^2 }
 \Big)^k. \]
Now set $k = \lfloor \frac{2V}{3(1+\eta)} \rfloor$,
and we obtain by Proposition \ref{Prop5}:
\[\Big( \frac{k \log\log T}{\eta^2 V^2 } \Big)^k 
   \le \exp\Big( -\frac23 V \log\frac{V}{\log\log T} 
   + \frac43 V\log\log V + \frac23 V \Big). \]
With
\begin{equation*}
(\log T)^2 (4 C)^k = \exp({\rm O}(V)),
 \text{ see \cite{RoBa},} 
\end{equation*} 
we get the assertion with an absolute ${\rm O}$-constant.
\qed

\Leer
\begin{proposition}
\label{Prop10}
(GRH) Let $\chi$ be a character mod $q$ with conductor $q_1$.
Further let $T$ be large, let
\[2(\log\log T)^{2} \le V \le \frac{\log T}{\log\log T},\] and let
$T\le t_1 <t_2 < \dots < t_R \le 2T$ be $V$-untypical ordinates
with $t_{r+1} - t_r \ge 1$ for all $1\le r< R$. Then
\[ R \ll T\exp\Bigl( -\frac23 V\log\frac{V}{\log\log T} + \frac43
V\log\log V + {\rm O}(V) \Bigr) \]
with an $O$-constant independent of $q$ and $\chi$.
\end{proposition}

\textbf{Proof:}

If $t$ is a $V$-untypical ordinate, 
then at least one of the criteria of Definition~\ref{Def1}
is false. For each criterion that is hurt, we give estimates for
the corresponding number $R_1$, $R_2$ and $R_3$ of such 
well-separated ordinates being counted in the Proposition. 

If criterion (i) is false for $t_r$, then there exists
a  $\sigma_r \ge \frac{1}{2}$ such that
\[ \Bigg|\underset{n\le x}{\sum}\frac{\Lambda(n)\chi_{1}(n)}{n^{\sigma_r +
      it_r} \log n}\frac{\log\frac{x}{n}}{\log x} \Bigg| > 2V,
\]
note here that $x =  T^{\frac{1}{V}}$.
  
The size of the sum over $n=p^{\alpha}$ with $\alpha \ge 2$ is
\begin{align*}
\Bigg|\underset{\substack{n=p^{\alpha}\le x \\ \alpha \ge
      2}}{\sum}\frac{\Lambda(n)\chi_{1}(n)}{n^{\sigma_r + it_r} \log
    n}\frac{\log\frac{x}{n}}{\log x} \Bigg|  
 &\le \underset{p\le \sqrt{x}}{\sum}\frac{1}{p} +
 \underset{\substack{p^{\alpha} \le x \\ \alpha \ge
     3}}{\sum}\frac{1}{p^{\frac{\alpha}{2}}} \\ 
 &\ll \log\log x \ll \log\log T \ll \sqrt{V}.
\end{align*}

So if we count the ordinates $t_{r}$ with 
\[ \Bigg| \underset{p\le x}{\sum}\frac{\chi_{1}(p)}{p^{\sigma_r +
      it_r}}\frac{\log\frac{x}{p}}{\log x}  \Bigg| \ge V,
\]
where again $x = T^{\frac{1}{V}}$, we get an upper bound for $R_{1}$.

Now we apply Proposition \ref{Prop4} of Maier and Montgomery,
we obtain
\begin{equation*}
   R_1V^{2k} \le \underset{r\le R}{\sum}\,
   \Bigg| \underset{p\le x}{\sum}
   \frac{\chi_{1}(p)}{p^{\sigma_r + it_r}}
   \frac{\log\frac{x}{p}}{\log x}  \Bigg|^{2k} 
   \ll T(\log T)^2 k!\, \Bigg(\underset{p\le x}{\sum}
   \frac{\log^2 \frac{x}{p}}{p\log^2 x}\Bigg)^k,
\end{equation*}
where $x^k \le T$ holds for every $k\le V$.

Now
\begin{equation*}
\underset{p\le x}{\sum}\frac{\log^2
    \frac{x}{p}}{p\log^2 x} \le \underset{p\le
  x}{\sum}\frac{1}{p} \ll \log\log x \le \log\log T .
\end{equation*}
As in \cite{RoBa}, 
we obtain with $k=\lfloor V \rfloor$:
\[ R_1 \ll T(\log T)^2 \Big(\frac{Ck \log\log T}{V^2}
\Big)^k = T\exp\Big(-V\log\frac{V}{\log\log T} +{\rm O}(V)
\Big). \]

\Leer
Now let (ii) be false, i.\,e.\ for $t_r$ there exists a $t_r'$ 
with  $|t_r - t_r'|\le 1$ and
\[ N\Big(t_r' + \frac{\pi\delta V}{\log(q_1T)},\chi\Big) -
N\Big(t_r' - \frac{\pi\delta V}{\log(q_1T)},\chi\Big) >
 (1+\delta) V. \]

With
\begin{equation*}
\delta V = \frac{\delta
  V}{\log(q_1T)}\log\left(\frac{q_1t_r'}{2\pi}\right) + {\rm o}(1) 
\text{ for }T\to\infty
\end{equation*}
we obtain
\[ N\Big(t_r' + \frac{\pi\delta V}{\log(q_1T)},\chi\Big) -
N\Big(t_r' - \frac{\pi\delta V}{\log(q_1T)},\chi\Big) -
\frac{\delta V}{\log(q_1T)}\log\Big(\frac{q_1t_r'}{2\pi}\Big)  \ge
V + {\rm O}(1).\]

Now we can apply Proposition \ref{Prop8}, if the $t_r'$ 
have a sufficiently large distance from another.
So instead of the sequence $t_r'$ being induced from $t_r$ for
$1\le r\le R_2$, consider the three subsequences
$t_{3s + \ell}'$ with $\ell\in\{1,2,3\}$,  $0\le s \le
\left\lfloor \frac{R_2-\ell}{3}\right\rfloor$, they
have the property $t_{3(s+1) + \ell}' - t_{3s +\ell}' \ge 1$. 
We can apply Proposition \ref{Prop8}
on any of the three subsequences and obtain
\[R_2 \le 3\Big(\Big\lfloor\frac{R_2}{3}\Big\rfloor +
  1\Big) + 2 \ll T\exp\Big(-\frac23 V\log\Big(\frac{V}{\log\log T}
  \Big) + \frac43 V\log\log V + {\rm O}(V) \Big).  \]

For $R_3$ we obtain, analogously as in  \cite{RoBa},
the same bound with a similar calculation.
\qed



\section{Logarithmic derivative of $L(s,\chi)$}
\label{sec6}

In this section, we consider only primitive characters. 

\Leer  
\begin{proposition}
\label{Prop11}
Let $\chi$ be a primitive character mod $q$,
$T$ be sufficiently large,
$\frac{1}{2}\le \sigma \le 2$, $T \le t \le2 T$
and $\Ls(\sigma +it,\chi)\not= 0$.
Then
\[\Re \frac{L'}{L}(\sigma +it,\chi) 
  = F(\sigma +it,\chi) - \frac{1}{2}\log(qT) + {\rm O}(1),\]
where $F(s,\chi) :=
\dis\underset{\rho}{\sum}\Re\frac{1}{s-\rho}$ and the sum runs through
all nontrivial zeros of $\Ls(s,\chi)$.
\end{proposition}

\textbf{Proof:}

We use the formula
\[\frac{L'}{L}(s,\chi) =
-\frac{1}{2}\log\frac{q}{\pi} - \frac{1}{2}\frac{\Gamma
  '}{\Gamma}\left( \frac{s + \mathsf{a}}{2} \right) + B(\chi) +
\underset{\rho}{\sum}\left( \frac{1}{s-\rho} + \frac{1}{\rho}
\right) \] 
that holds for primitive characters,
where $\Re B(\chi) = - \underset{\rho}{\sum}\Re(\frac{1}{\rho})$ and the sum
runs through all nontrivial zeros $\rho$ of $\Ls(s,\chi)$. By
Stirling's formula we obtain
\begin{align*}
\Re\frac{L'}{L}(\sigma +it,\chi) 
  &= -\frac{1}{2}\log\frac{q}{\pi} -
  \frac{1}{2}\Re\frac{\Gamma '}{\Gamma}\left( \frac{\sigma +it +
   \mathsf{a}}{2} \right) + \Re B(\chi)+
   \underset{\rho}{\sum}\Re\left( \frac{1}{\sigma +it - \rho} +
   \frac{1}{\rho} \right) \\  
  &= -\frac{1}{2}\log q - \frac{1}{2}\log|\sigma +
    it + \mathsf{a}| + F(\sigma + it, \chi) +
    {\rm O}(|\sigma + it + \mathsf{a}|^{-1}) +
    {\rm O}(1) \\  
  &= F(\sigma + it, \chi) -
    \frac{1}{2}\log(qT) + {\rm O}(1).
\end{align*}
\qed

\Leer
\begin{proposition}
\label{Prop12}
Let $\chi$ be a primitive character mod $q$.
Let $x\ge 1$, and consider $z\in\C$ that is not a pole of 
$\frac{L'}{L}(z,\chi)$. Then
\[\underset{n\le
  x}{\sum}\frac{\chi(n)\Lambda(n)}{n^{z}}\log\Big(\frac{x}{n}\Big) = -
\frac{L'}{L}(z,\chi) \log x - \Bigl( \frac{L'}{L} \Bigr) '
(z,\chi) - \underset{\rho}{\sum}\frac{x^{\rho - z}}{(\rho - z)^{2}} -
\underset{n \geq 0}{\sum}\frac{x^{-2n -\mathsf{a} - z}}{(z + 2n +
  \mathsf{a})^{2}}.\]
\end{proposition}

\textbf{Proof:}
Since
\[\frac{\Ls '}{\Ls}(s,\chi) \ll \log(q|s|) ~\text{for}~ \Re s \le
-\frac{1}{2} \text{ and } |s + m| > \frac{1}{4} \text{ for all } m\in\N,\]
the proof works analogously to \cite{RoBa}, 
where the term coming from the pole at $s=1$ is removed
and the sum over the trivial zeros has been adjusted.
\qed

\Leer
Estimating the last sum analogously to \cite{RoBa}, we obtain:

\begin{proposition}
\label{Prop13}
Let $\chi$ be a primitive character mod $q$, 
$T\ge 1$ and  $1\le x\le T$. 
Let $z\in\C$, $\Re z\ge 0$, $T\le\Im z\le2T$, and let
$z$ be not a pole of $\frac{\Ls'}{\Ls}(z,\chi)$.

Then
\begin{equation}
\label{eProp13}
\underset{n\le
   x}{\sum}\frac{\chi(n)\Lambda(n)}{n^{z}}\log\Big(\frac{x}{n}\Big) = -
   \frac{L'}{L}(z,\chi) \log x - \Bigl( \frac{L'}{L} \Bigr) '
   (z,\chi) - 		\underset{\rho}{\sum}\frac{x^{\rho -
   z}}{(\rho - z)^{2}} + {\rm O}(T^{-1}).
\end{equation}
\end{proposition}
\qed


\section{Lower bound for $\log|L(s,\chi)|$}
\label{sec7}

With the aid of $V$-typical ordinates, we estimate
$\log\Ls(s,\chi)$ from below.
 
\begin{proposition}
\label{Prop14}
(GRH) Let $\chi$ be a nonprincipal character mod $q$
induced by $\chi_{1}$ mod $q_1$. 
Let $T$ be sufficiently large and $T\le t \le2T$.

Then for all $\frac{1}{2}\le \sigma \le 2$ and $2\le x \le T$ it holds
that
\[\log |L(\sigma + it, \chi)| \ge \Re\Big(\underset{n\le x}{\sum}
\frac{\Lambda(n) \chi_{1}(n)}{n^{\sigma + it}\log
  n}\frac{\log\frac{x}{n}}{\log x}\Big) - \Big( 1 +
\frac{x^{\frac{1}{2} - \sigma}}{(\sigma - \frac{1}{2})\log
  x}\Big)\frac{F(\sigma + it,\chi)}{\log x} +
{\rm O}\Big(\sqrt{\frac{\log q}{\log\log q}}\Big),\]
where $F$ is the function from Proposition \ref{Prop11}.
\end{proposition}

\textbf{Proof:}
At first, let $\chi$ be primitive.
By integrating equation (\ref{eProp13}) from $z = \sigma +it$
to $z = 2 +it$, we obtain analogously to \cite{RoBa}: 
\[\log |L(\sigma + it, \chi)| \ge
 \Re\Big(\underset{n\le x}{\sum}
\frac{\Lambda(n) \chi(n)}{n^{\sigma + it}\log
  n}\frac{\log\frac{x}{n}}{\log x}\Big) - \Big( 1 +
\frac{x^{\frac{1}{2} - \sigma}}{(\sigma - \frac{1}{2})\log
  x}\Big)\frac{F(\sigma + it,\chi)}{\log x} + {\rm O}(1).
\] 

Now let $\chi$ mod $q$ be not primitive and
induced by the primitive character $\chi_1$ mod $q_1$ .

Then we have 
\begin{equation}
\label{primitiv} 
  \Ls(s,\chi) = L(s,\chi_{1}) 
  \prod_{p\mid q}\Big(1-\frac{\chi_1(p)}{p^s}\Big).
\end{equation}

We obtain with equation (\ref{primitiv}): 
\begin{align*}
\log\Big|\Ls(s,\chi)\Big| &= \log|\Ls(s,\chi_1)| +
\underset{p|q}{\sum}\log\left|1-\frac{\chi_{1}(p)}{p^s}\right| \\
&\geq  \Re\Big(\underset{n\le x}{\sum}
\frac{\Lambda(n) \chi_{1}(n)}{n^{\sigma + it}\log
  n}\frac{\log\frac{x}{n}}{\log x}\Big) - \Big( 1 +
\frac{x^{\frac{1}{2} - \sigma}}{(\sigma - \frac{1}{2})\log
  x}\Big)\frac{F(\sigma + it,\chi_{1})}{\log x} + {\rm O}(1)
+\underset{p|q}{\sum}\log\left|1-\frac{\chi(p)}{p^s}\right|.
\end{align*}

For the last sum we get
\begin{equation}
\label{nichtprimitiv}
\underset{p|q}{\sum}\log\left|1-\frac{\chi(p)}{p^s}\right| \le
\underset{p|q}{\sum}\frac{1}{p^{1/2}} \le
\underset{j=1}{\overset{2\log q}{\sum}}\frac{1}{p_j^{1/2}} \ll
\sqrt{\frac{\log q}{\log\log q}}.
\end{equation}
From equation (\ref{primitiv}) we see further that
\[F(s,\chi) = F(s,\chi_1),\]
so we get the stated bound.
\qed

\Leer
Now we would like to give an estimate for
$\Ls(s,\chi)$ in the interval $\Re(s)\in\left(\frac12 , 2\right)$. 
For this, we split the interval at 
$\frac12 + \frac{V}{\log T}$ and give a bound for 
each part. This is done in the next two propositions.

\begin{proposition}
\label{Prop15}
(GRH) Let $\chi$ be a nonprincipal
character mod $q$, and further let $T$ be sufficiently large, at least $T\ge q$,
let $V\in\lbrack (\log\log T)^{2}, \frac{\log T}{\log\log T}\rbrack$
and let $t \in \lbrack T , 2T \rbrack$ be $V_{\delta,\chi,q}$-typical of
order $T$.

Then it holds for $\frac{1}{2} + \frac{V}{\log T}\le\sigma\le 2$, that
\[\log |L(\sigma + it,\chi)| \ge f_{\delta, q}(V,\sigma +it),\]
where $f_{\delta, q}: \R\times\C \rightarrow \R$,
$f_{\delta, q}(V,\sigma +it) = {\rm O}\left(\frac{V}{\delta} +
  \sqrt{\frac{\log q}{\log\log q}}\right)$. 
\end{proposition}

\textbf{Proof:}

In Proposition \ref{Prop14} 
we set $x=T^{\frac{1}{V}}$. Then $2\le x \le T$, and since 
$\frac{1}{2} + \frac{V}{\log T}\le\sigma$, we have
\[\frac{x^{\frac{1}{2} - \sigma}}{(\sigma - \frac{1}{2})\log x} \le
\frac{\exp(-V\frac{\log x}{\log T})}{V\frac{\log x}{\log T}}  = e^{-1}
\le 1.\]
Applying now Proposition \ref{Prop14}, we obtain:
\begin{equation*}
\log\left|L(\sigma + it, \chi)\right| \ge -2V -2\frac{V}{\log
  T}F(\sigma +it,\chi) + {\rm O}\Big(\sqrt{\frac{\log q}{\log\log
      q}}\Big) =: f_{\delta, q}(V,\sigma +it), 
\end{equation*}
since $t$ is $V$-typical.

We aim to majorize $F(\sigma +it,\chi)$ independent from $q$ and
$\chi$. As in \cite{RoBa}, we divide the region of the zero-ordinates
in two parts as follows.
\begin{itemize}
  \item[(i)] $\gamma$ with $\frac{2\pi n \delta
      V}{\log(q_1T)}\le |t-\gamma|\le \frac{2\pi (n+1) \delta
      V}{\log(q_1T)} \text{ for } 0\le n \le N = \left\lfloor 
	\frac{\log(q_1T)}{4\pi \delta V}\right\rfloor,$
\item[(ii)] 
       $\gamma$ with $\{ \gamma ~:~ |\gamma - t|\ge \frac{1}{2} \}$,
     where $q_{1}$ denotes the conductor of $\chi$ mod $q$.
\end{itemize}

Consider the set of $\gamma$ from (i):
\begin{align*}
\underset{\gamma\text{ from (i)}}{\sum}\Re \frac{1}{\sigma +it -\frac12
  - i\gamma} &=2\underset{\gamma\text{ from (i)}}{\sum}\frac{(\sigma
  -\frac12)}{(\sigma -\frac12)^2 + (t-\gamma)^2}\\ 
&\le 2(1+\delta )V\underset{n=0}{\overset{N}{\sum}}\frac{(\sigma
  -\frac{1}{2})}{(\sigma -\frac{1}{2})^{2} +
  \big(\frac{2\pi n \delta V}{\log(q_1T)}\big)^{2}} \text{ since } t
  \text{ is } V\text{-typical, (ii),}\\ 
&\le 4V\Bigl(\frac{1}{\sigma - \frac{1}{2}} +
\frac{\log(q_1T)}{4\delta V} \Bigr), 
\end{align*}
since for $a,c\in\R_{>0}$ and $N\in\N$ we have
 $\underset{n=0}{\overset{N}{\sum}} \frac{a}{a^2 +(cn)^2}
 \le \frac1a + \frac{\pi}{2c} $, see \cite{RoBa} Prop. 6, and
we continue with
\begin{align*} 
&\le 4\log(q_1T) + \frac{\log (q_1T)}{\delta }
 \le 5\frac{\log(qT)}{\delta}.
\end{align*}

For the sum over $\gamma$ with (ii) we work with
the known formula 
\begin{equation}\label{s.86}
\underset{\rho\in\mathcal{N}(\chi)}{\sum}
   \frac{1}{1+(t-\Im(\rho))^2}\ll\log(q(2+|t|))  
\end{equation}
holding for primitive characters mod $q$.
Since $\mathcal{N}(\chi) = \mathcal{N}(\chi_1)$ if $\chi$ mod $q$
is induced by $\chi_1$ mod $q_1 \le q$, we can use this formula also
in the case of a nonprimitive character mod $q$. 

For $0\le\sigma-\frac12\le\frac32$ and $|t-\gamma|\ge\frac12$ we have
\begin{equation}\label{ungl}
\frac{\sigma - \frac{1}{2}}{(\sigma - \frac{1}{2})^{2}
   + (t - \gamma )^{2}} \le \frac{8}{1+(t-\gamma)^2},
\end{equation}
therefore we can estimate the sum over $\gamma$ with (ii) 
using (\ref{ungl}) and (\ref{s.86}) as follows:
\begin{align*}
\underset{|t-\gamma|\ge \frac12}{\sum}\Re\Big(\frac{1}{\sigma + it -
    \frac12 - i\gamma}\Big) &=\underset{|t-\gamma|\ge
  \frac12}{\sum}\frac{\sigma-\frac12}{(\sigma-\frac12)^2 +
  (t-\gamma)^2} \\ 
&\le\underset{|t-\gamma|\ge
  \frac12}{\sum}\frac{8}{1+(t-\gamma)^2} 
\le\underset{\rho\in\mathcal{N}(\chi)}{\sum}
  \frac{8}{1 + (t-\Im(\rho))^2}
\ll \log(qt).
\end{align*}

Now consider $g(x):= \frac{\log(qx)}{\log x}$, we see that
$g(x)$ is monotonously decreasing for $x>1$, and so for
$x\ge q$ we have $g(x)\le g(q) = 2$.

We resume the two results for the regions (i) and (ii) as follows:
\begin{align*}
  &\Big|2\frac{V}{\log T}F(\sigma +it,\chi)\Big| \ll
     \frac{\log(qT)}{\log T}\frac{V}{\delta} \ll \frac{V}{\delta}
     \text{ since }q\leq T,
\end{align*}
which gives the asserted bound for $f_{\delta,q}(V,\sigma+it)$.
\qed

\Leer
\begin{proposition}
\label{Prop16}
(GRH) Let $\chi$ be a character mod $q$, let $T$ be sufficiently large,
$V \in\lbrack (\log\log T)^{2}, \frac{\log T}{\log\log T}\rbrack$
and  $t \in \lbrack T , 2T \rbrack$ be $V$-typical (of order $T$).

Then we have for all $\frac{1}{2} < \sigma \le \sigma_{0} 
  = \frac{1}{2} + \frac{V}{\log T}$:
\[ \log |L(\sigma + it,\chi)| \ge
 \log |L(\sigma_{0} + it,\chi)|
  - V\log\frac{\sigma_{0}-\frac{1}{2}}{\sigma - \frac{1}{2}} 
   - 2(1 + \delta)V\log\log V 
   + {\rm O}\Big(\frac{V}{\delta^{2}}
    + \sqrt{\frac{\log q}{\log\log q}}\Big).\]
\end{proposition}

\textbf{Proof:}

Consider at first a primitive character $\chi$ mod $q$, i.\,e.\ $q_1 = q$.
We work as in \cite{RoBa}, p.\ 8, and get:
\begin{equation*}
\log|\Ls (\sigma_{0} +it, \chi)| - \log|\Ls (\sigma + it,\chi)|
   \le \frac{1}{2}\underset{\gamma}{\sum}
   \log\frac{(\sigma_{0} - \frac{1}{2})^{2}
    + (t - \gamma)^{2}}{(\sigma - \frac{1}{2})^{2} 
   + (t - \gamma)^{2}}.
\end{equation*}

In order to estimate the sum, we divide the set of $\gamma$ in three
subsets such that we can make use of the fact that $t$ is a
$V$-typical ordinate.

The division of the $\gamma$ is as follows.
\begin{enumerate}
\item[(a)] $\gamma$ with $|t-\gamma | \le \frac{\pi V}{\log V
    \log(qT)}$,
\item[(b)] $\gamma$ with $\big(2\pi \delta n +
   \frac{\pi}{\log V} \big) \frac{V}{\log (qT)} 
   \le |t-\gamma | \le \big(2\pi \delta (n+1)
    + \frac{\pi}{\log V} \big) \frac{V}{\log (qT)} 
   \quad\big(0\le n \le N = \big\lfloor \frac{\log (qT)}{4\pi \delta V}
      \big\rfloor\big)$,
\item[(c)] $\gamma$ with $\{ \gamma \;:\; |t-\gamma | \ge\frac{1}{2} \}$.
\end{enumerate}

Since $\sigma\le\sigma_0$, we have
\[\frac{(\sigma_{0} - \frac{1}{2})^{2} + (t - \gamma)^{2}}{(\sigma -
  \frac{1}{2})^{2} + (t - \gamma)^{2}} \le \frac{(\sigma_{0} -
  \frac{1}{2})^{2}}{(\sigma - \frac{1}{2})^{2}}. \]
 
For the $\gamma$ from (a) we use property (iii) from the definition of
$V$-typical and obtain
\begin{equation*}
  \frac12 \underset{|t-\gamma | \le
   \frac{\pi V}{\log V \log(qT)}}{\sum}
   \log\frac{(\sigma_{0} - \frac{1}{2})^{2} 
    + (t - \gamma)^{2}}{(\sigma - \frac{1}{2})^{2} 
    + (t - \gamma)^{2}} \le \frac12 
    \underset{|t-\gamma | \le \frac{\pi V}{\log V
        \log(qT)}}{\sum}\log\frac{(\sigma_{0} -
        \frac{1}{2})^{2}}{(\sigma - \frac{1}{2})^{2}} 
     \le V\log\frac{\sigma_{0} -\frac{1}{2}}
    {\sigma - \frac{1}{2}}.
\end{equation*}

We use the fact that $\frac{(\sigma_{0} - \frac{1}{2})^{2}
  + (t - \gamma)^{2}}{(\sigma - \frac{1}{2})^{2} + (t - \gamma)^{2}}$
is decreasing in $|t-\gamma|$. With this, we estimate the set of $\gamma$
in (b) using property (ii) in the definition of $V$-typical.
For the $\gamma$ with (c) we use the general zero estimate for $\Ls(s,\chi)$
and obtain in the same way as in \cite{RoBa}: 
\begin{equation*}
  \frac{1}{2}\underset{\gamma\text{'s in (b)}}{\sum}
  \log\frac{(\sigma_{0} - \frac{1}{2})^{2} + (t -
    \gamma)^{2}}{(\sigma - \frac{1}{2})^{2} + (t - \gamma)^{2}}
  \le 2(1+\delta ) V \log\log V + {\rm O}\left(\frac{V}{\delta^{2}}
  \right)  
\end{equation*}
and
\begin{equation*}
\frac{1}{2}\underset{|\gamma -
  t|\ge\frac{1}{2}}{\sum}\log \frac{(\sigma_{0} -
    \frac{1}{2})^{2} + (t - \gamma)^{2}}{(\sigma - \frac{1}{2})^{2} +
    (t - \gamma)^{2}} \ll \frac{V}{\log\log T}. 
\end{equation*}
 
This gives the assertion for primitive characters.

Now if $\chi$ is not primitive mod $q$ and induced by the primitive
character $\chi_1$ mod $q_{1}$, we use equation (\ref{nichtprimitiv})
and obtain
\begin{align*}
\log|\Ls(\sigma +it,\chi)| &= \log|\Ls(\sigma +it,\chi_1)| +
{{\rm O}\Big(\sqrt{\frac{\log q}{\log\log q}}\Big)}\\ 
  &\ge \log|\Ls(\sigma_0 +it,\chi_1)|
- V\log \frac{\sigma_{0}-\frac12}{\sigma-\frac12} 
-2(1+\delta)V\log\log V+ {\rm O}(\frac{V}{\delta^{2}})
+{\rm O}\Big(\sqrt{\frac{\log q}{\log\log q}}\Big) \\ 
&= \log|\Ls(\sigma_0 +it,\chi))| +
- V\log \frac{\sigma_{0}-\frac12}{\sigma-\frac12} 
-2(1+\delta)V\log\log V+ {\rm O}\Big(\frac{V}{\delta^{2}}
+\sqrt{\frac{\log q}{\log\log q}}\Big).
\end{align*}
\qed

\Leer
At the end of this section we combine the results from propositions
\ref{Prop9}, \ref{Prop15} and \ref{Prop16}. 
With these, we obtain a lower bound for the whole
stripe $\Re(s)\in\left(\frac12,2\right)$.

\begin{proposition}
\label{Prop17}
(GRH) Let $\chi$ be a character mod $q$, $|t|$ be sufficiently large,
at least $|t|\ge q$, and $\frac{1}{2}< \sigma \le 2$. Then
\[ \log |\Ls(\sigma + it, \chi)| \ge - \frac{\log |t|}{\log\log
  |t|}\log\frac{1}{(\sigma - \frac{1}{2})} - 3 \frac{\log |t|
  \log\log\log |t|}{\log\log |t|}. \]
\end{proposition} 

\textbf{Proof:} 

As in \cite{RoBa}, we choose
\[V= \frac{\log|t|}{\log\log|t|} \text{ and } \delta = \frac{1}{2}, \]
note that then ${\rm O}\left(\frac{V}{\delta^2} +
  \sqrt{\frac{\log q}{\log\log q}}\right) = {\rm O}\left( V
\right)$. 
\qed

\Leer
By now, we gave estimates for $\Ls(s,\chi)$ in a region for
sufficiently large $\Im(s)$. We also need an estimate for
$\Ls(s,\chi)$ in the remaining region, which we give in the next
Proposition.

\begin{proposition}
\label{Prop18}
(GRH) Let $x$ be large, $c>0$. Further let 
$T_0(x) := T_0 := 2^{\left\lfloor (\log x)^{3/5}(\log\log
    x)^c\right\rfloor}$, and $\sigma = \frac12 + \frac{1}{\log x}$.
Then there exists a $C>0$, such that for all 
$|t|\le T_0$, $q\le \sqrt{T_0}$ and a nonprincipal character 
$\chi$ mod $q$ we have
\[\left|\Ls(\sigma +it,\chi)\right|\ge T_0^{-C\log\log x}. \]
\end{proposition}

\textbf{Proof:}

At first, let $\chi$ be a primitive character mod $q$, 
and $q\le\sqrt{T_0}$.
By the explicit formula for the logarithmic derivation 
of $L$ we obtain
\begin{equation*}
\underset{\sigma
  +it}{\overset{2+it}{\int}}\frac{\Ls'}{\Ls}(s+it,\chi)ds =
\underset{\sigma
  +it}{\overset{2+it}{\int}} \Big(
   \underset{\substack{\rho\in\mathcal{N}(\chi)
    \\ |\Im(s)- \Im(\rho)|\le 1}}{\sum}\frac{1}{s-\rho} +
{\rm O}(\log(q(2+|\Im(s)|))) \Big) ds,
\end{equation*}
hence
\begin{align*}
 &\log\Ls(2+it,\chi) - \log\Ls(\sigma+it,\chi) \\
&=\underset{\substack{\rho\in\mathcal{N}(\chi) \\ |t- \Im(\rho)|\le
    1}}{\sum}\log(2+it-\rho) -
\underset{\substack{\rho\in\mathcal{N}(\chi) \\ |t- \Im(\rho)|\le
    1}}{\sum}\log(\sigma +it -\rho) + {\rm O}(\log(q(2+|t|))). 
\end{align*}
Considering the real parts, it follows that
\[\log|\Ls(\sigma+it,\chi)|^{-1} =
\underset{\substack{\rho\in\mathcal{N}(\chi) \\ |t-
      \Im(\rho)|\le 1}}{\sum} \log\left|\frac32
    +i(t-\Im(\rho))\right| +
\underset{\substack{\rho\in\mathcal{N}(\chi) \\ |t-
      \Im(\rho)|\le 1}}{\sum}\log\frac{1}{|\sigma+it-\rho|} +
{\rm O}(\log(q(2+|t|))). \]
To give an estimate of the first sum,
we have
\[\left|\frac32 + i(t-\Im(\rho))\right| 
  \le\frac52 \text{ for } |t-\Im(\rho)|\le 1,
\text{ hence } \underset{|t-\Im(\rho)|\le
  1}{\sum}\log\left|\frac32+i(t-\Im(\rho))\right| \ll \log(qt), \]
and to give an estimate offor the second sum, we have
\begin{align*}
&|\sigma +it -\rho|^{-1} = \left|\frac{1}{\log x} +
  i(t-\Im(\rho))\right|^{-1} \le \log x, \\ 
&\text{ hence } \underset{|t-\Im(\rho)|\le
  1}{\sum}\log\frac{1}{|\sigma+it-\rho|} \ll \log(qt) \log\log x.
\end{align*}

Therefore we obtain
\[\log|\Ls(\sigma+it,\chi)|^{-1} \ll \log(qt) \log\log x.  \]

If we note that $t\le T_0$ and $q\le \sqrt{T_0}$, we obtain
\[\log|\Ls(\sigma+it,\chi)|^{-1} \ll \log T_0 \log\log x. \]
This gives the assertion for primitive characters.

Now let $\chi$ be a nonprimitive character mod $q$ and induced by 
$\chi_1$ mod $q_1$. We conclude:
\begin{align*}
\log|\Ls(\sigma+it,\chi)|^{-1} &= \log|\Ls(\sigma+it,\chi_1)|^{-1} -
\sum_{p|q}\log\Big|1-\frac{\chi(p)}{p^s}\Big| \\ 
	&=\log|\Ls(\sigma+it,\chi_1)|^{-1} +
        {\rm O}\Big(\sqrt{\frac{\log T_0}{\log\log T_0}}\Big) \\ 
	&\ll \log T_0 \log\log x
        \Big(1+{\rm O}\Big(\frac{1}{\sqrt{\log T_0 \log\log
                T_0}\log\log x}\Big)\Big)\\ 
	&\ll \log T_0 \log\log x.
\end{align*}
\qed


\section{Majorant of $\left|x^{z}\Ls (z,\chi)^{-1}\right|$}
\label{sec8}

In this section we give a majorant of
$\left|x^{z}\Ls (z,\chi)^{-1}\right|$
for certain $z$.
It is a consequence of Propositions \ref{Prop15} and \ref{Prop16}.

\begin{proposition}
\label{Prop19}
(GRH) Let $\chi$ be a character mod $q$. Further let
$t$ be sufficiently large (at least $t\ge q$),
$x\ge t $, $V' \in \left\lbrack (\log\log t)^{2}, 
   \frac{\log (t/2)}{\log\log (t/2)} \right\rbrack$,
$V\ge V'$,
$t$ be $V'$-typical of order $T'$.

Then for $V'\le (\Re z -\frac{1}{2})\log x \le V$, $|\Im z| = t$, we
have
\[\Big|x^{z}\Ls(z,\chi)^{-1}\Big| \le \sqrt{x} \exp\Big(
  V\log\frac{\log x}{\log t} + 2(1+\delta )V\log\log V  +
  {\rm O}\Big(V\delta^{-2} + \sqrt{\frac{\log x}{\log\log
        x}}\Big)\Big).\]
\end{proposition}

\textit{Proof:}
By taking notion of the changed error term,
everything remains as in \cite{RoBa}, see Proposition 22
there.
\qed


\section{Upper bound for $M(x,q,a)$}
\label{sec9}

We need some preliminaries for the proof of the theorem. 

For a character $\chi$ mod $q$, let

\[A(x,\chi,q) := \frac{1}{2\pi i}\underset{1 + \frac{1}{\log x} - i
  2^K}{\overset{1 +\frac{1}{\log x} + i 2^K}{\int}}
  \frac{x^s}{\Ls(s,\chi) s}ds,
\text{ where } K:=\left\lbrack\frac{\log x}{\log 2}\right\rbrack,
 \] 
and by Perron's formula we have:
\begin{equation}
\label{vorueb}
  M(x,q,a) = \frac{1}{\varphi(q)}\underset{\chi (q)}
    {\sum}\overline{\chi}(a)A(x,\chi, q) + {\rm O}(\log x).
\end{equation}

We aim to give a good upper bound for $A(x,\chi,q)$.

Further we assume w.l.o.g., that
$x\ge q^2$, as otherwise we can estimate trivially.

Now we give some definitions being valid during this section.

\begin{definition}
\begin{align*}
K&:=\left\lbrack\frac{\log x}{\log 2}\right\rbrack, 
\kappa:=\left\lfloor (\log x)^{3/5} (\log\log x)^{c}\right\rfloor,\\
T_{k}&:=2^{k} \text{ for }
  \kappa\leq k\leq K, \text{ so } q^2\le T_{\kappa}\le T_k.
\end{align*}
For $k$ with $\kappa\le k < K$ and for
   $n \in \N\cap \lbrack T_k, 2T_k)$, we define the integer $V_n$
   to be the smallest integer in the interval
   $\left\lbrack(\log\log T_k)^2 +1 , \frac{\log T_k}{\log\log T_k}
   \right\rbrack$, such that all points in 
   $\lbrack n, n+1\rbrack$ are $V_n$-typical ordinates of order $T_k$.
The existence of these $V_n$  is obtained by Proposition \ref{Prop9}.
\end{definition}

\begin{lemma}
\label{l1}
Let $x\ge 2 \text{, } c>1$, $q\in\N$ and 
$1<q\le 2^{\kappa/2}$.
Further let $\chi$ be a nonprincipal character mod $q$ and
$\delta\in\left(0,1\right\rbrack$. Then
\[\frac{A(x,\chi,q)}{\sqrt{x}} \ll_{\delta} \exp\Big((\log
    x)^{3/5}(\log\log x)^{c + 1 +
    \delta}\Big) + B(x,\chi,q), \]
where
\begin{equation*}
B(x,\chi,q) = \underset{n = T_{\kappa}}{\overset{T_K - 1}{\sum}}
  \frac{1}{n}
\exp\Big(V_n \log\Big(\frac{\log x}{\log
      n}\Big) + 2(1 + 2\delta)V_n\log\log V_n + D \sqrt{\frac{\log
      x}{\log\log x}} \Big)
\end{equation*} 
with an absolute constant $D>0$. 
\end{lemma}

\textbf{Proof:}

We choose the following path of integration $S(x,\chi,q)$,
we describe it for the upper half plane
$\Im(z)\ge 0$, it passes out analogously in the lower half plane.

\begin{enumerate}
\item A vertical segment
 $\dis \left\lbrack \frac{1}{2} + \frac{1}{\log x},
 \frac{1}{2} + \frac{1}{\log x} + iT_{\kappa} \right\rbrack.$ 
\item  Further vertical segments
   $\dis \left\lbrack \frac{1}{2} + \frac{V_n}{\log x} + in, \frac{1}{2}
    + \frac{V_n}{\log x} +i(n+1)\right\rbrack$.
\item A horizontal segment
$\dis \left\lbrack \frac{1}{2} + \frac{1}{\log x} + iT_{\kappa}
  ,\frac{1}{2} + \frac{V_{T_{\kappa}}}{\log x} + iT_{\kappa}
\right\rbrack$.
\item Additional horizontal segments for $T_{\kappa} \le n \le T_K -
  2$, namely

$\dis \left\lbrack \frac{1}{2} + \frac{V_n}{\log x} + i(n + 1)
  ,\frac{1}{2} + \frac{V_{n + 1}}{\log x} + i(n + 1) \right\rbrack$.
\item The last horizontal segment
$\dis \left\lbrack \frac{1}{2} + \frac{V_{T_K -1}}{\log x} + iT_K  ,1
  + \frac{1}{\log x} + iT_K \right\rbrack$.
\end{enumerate}
Hence
\[ |A(x,\chi, q)| = \frac{1}{2\pi}\Big|\underset{S(x,\chi,
    q)}{\int}\frac{x^s}{\Ls(s,\chi) s}ds \Big|. \]

We consider just the first segment more accurate,
the others can be estimated analogously to
\cite{RoBa}:
 
Ad 1.:
\begin{align*}
\frac{1}{2\pi}\Biggl|\underset{\substack{S(x,\chi, q)\\ |\Im(z)|\le
  T_{\kappa}}}{\int}\frac{x^s}{\Ls(s,\chi) s}ds \Biggr|  
  &\le \frac{1}{2\pi}x^{\frac12 +\frac{1}{\log
   x}}\underset{-T_{\kappa}}{\overset{T_{\kappa}}{\int}}\Big|\Ls\Big(\frac{1}{2}+
   \frac{1}{\log x} + it,
    \chi\Big)\Big|^{-1}\frac{dt}{\sqrt{\frac{1}{4} + t^2}} \\ 
   &\le \frac{e}{2\pi} \sqrt{x} \underset{|t| \le
    T_{\kappa}}{\max}\Big|\Ls\Big(\frac{1}{2}+ \frac{1}{\log x}
    + it, \chi\Big)\Big|^{-1}
   \underset{-T_{\kappa}}{\overset{T_{\kappa}}{\int}}\frac{dt}{\sqrt{\frac{1}{4}
    + t^2}} \\ 
   &\le \sqrt{x} \underset{|t| \le
   T_{\kappa}}{\max}\Big|\Ls\Big(\frac{1}{2}+ \frac{1}{\log x}
   + it, \chi\Big)\Big|^{-1}
  \underset{0}{\overset{T_{\kappa}}{\int}}\frac{dt}{\sqrt{\frac{1}{4}
  + t^2}} \\ 
  &\le 2\sqrt{x} \underset{|t| \le
    T_{\kappa}}{\max}\Big|\Ls\Big(\frac{1}{2}+ \frac{1}{\log x}
     + it, \chi\Big)\Big|^{-1} \log T_{\kappa}\\ 
  &\ll\sqrt{x} \left(\log T_{\kappa}\right)  
      T_{\kappa}^{C\log\log x} \text{ by Prop. }\ref{Prop18} \\ 
  &\le \sqrt{x}T_{\kappa}^{C_1\log\log x} \text{ with }C_{1}=C+1.
\end{align*}

Ad 2.:\\
\begin{align*}
  \frac{1}{2\pi} \Biggl|\underset{\frac{1}{2} + \frac{V_n}{\log x} +
    in}{\overset{\frac{1}{2} + \frac{V_n}{\log x}
      +i(n+1)}{\int}}&\frac{x^s}{\Ls(s,\chi) s}ds \Biggr|
  \le \frac{1}{2\pi n}
  \max_{\substack{z\in\{\frac{1}{2} + \frac{V_n}{\log x} + it;  \\ 
    t\in\lbrack n, n+1\rbrack\}  }}
  \Big|x^z\Ls(z,\chi)^{-1}\Big| 
   \text{ as } |s|\ge|n|     \\ 
  &\le \frac{1}{n}\sqrt{x}\exp\Big(V_n \log\Big(\frac{\log x}{\log
        n}\Big) + 2(1 + \delta)V_n\log\log V_n +
    D\Big(\frac{V_n}{\delta^2} + \sqrt{\frac{\log x}{\log\log
          x}}\Big) \Big), 
\end{align*}
where $D>0$ is an absolute constant,
see Proposition \ref{Prop19}.

Ad 3.:\\
\begin{equation*}
\frac{1}{2\pi} \Biggl|\underset{\frac{1}{2} + \frac{1}{\log x} +
    iT_{\kappa}}{\overset{\frac{1}{2} + \frac{V_{T_{\kappa}}}{\log x}
      + iT_{\kappa}}{\int}}\frac{x^s}{\Ls(s,\chi) s}ds \Biggr| \le
\sqrt{x}T_{\kappa}^3 \text{ by Prop. }\ref{Prop17}.
\end{equation*}

Ad 4.:\\
Here we use Proposition \ref{Prop19} for $n$ with
$T_{\kappa} \le n \le T_K -2 $:
\begin{align*}
\Biggl|&\underset{\frac{1}{2} + \frac{V_n}{\log x} + i(n +
    1)}{\overset{\frac{1}{2} + \frac{V_{n + 1}}{\log x} + i(n +
      1)}{\int}}\frac{x^s}{\Ls(s,\chi) s}ds \Biggr| \\
 &\le \frac{1}{n}\sqrt{x}\exp\Big(V_n \log\Big(\frac{\log x}{\log
       n}\Big) + 2(1 + \delta)V_n\log\log V_n +
   D\Big(\frac{V_n}{\delta^2} + \sqrt{\frac{\log x}{\log\log
         x}}\Big) \Big) \\ 
&+\frac{1}{n +1}\sqrt{x}\exp\Big(V_{n+1} \log\Big(\frac{\log
      x}{\log(n+1)}\Big) + 2(1 + \delta)V_{n+1}\log\log V_{n+1} +
  D\Big(\frac{V_{n+1}}{\delta^2} + \sqrt{\frac{\log x}{\log\log
        x}}\Big) \Big).
\end{align*}

Ad 5.:\\
We obtain using Proposition \ref{Prop15}:
\begin{equation*}
\frac{1}{2\pi} \Biggl|\underset{\frac{1}{2} + \frac{V_{T_K -1}}{\log x}
    + iT_K}{\overset{1 + \frac{1}{\log x} +
      iT_K}{\int}}\frac{x^s}{\Ls(s,\chi) s}ds \Biggr| 
	\le_{\delta} \sqrt{x}.
\end{equation*}
\qed

\Leer
The following proposition is similar to Proposition 23 in
\cite{RoBa}, the modification here is necessary, but the
proof works analogously.

\begin{proposition}
\label{Prop20}
Let $A,C>0$ and let $A\ge 4C^4 +1$, then for $V>e^{3C/2}$
it holds that
\[AV - \frac23 V \log V + CV\log\log V \le e^{3A/2} \Big(\frac32
  A\Big)^{3C/2}.\]
\end{proposition}

\begin{lemma}
\label{l2}
Under the conditions of Lemma \ref{l1} we have
\[B(x,\chi,q) \ll_{\delta} \exp\Big((\log x)^{3/5}(\log\log
  x)^{13/2 - 3c/2 + 8\delta} \Big).  \]
\end{lemma}

\textbf{Proof:}

We define for $\kappa \le k < K$:
\[ B(T_k,x,\chi,q) := \underset{T_k \le n <
  2T_k}{\sum}\frac{1}{n}\exp\Big(V_n\log\Big(\frac{\log x}{\log
      n}\Big) + 2(1+2\delta)V_n\log\log V_n\Big), \]
then
\[ B(x,\chi,q) \le K \underset{\kappa\le k <
  K}{\max}B(T_k,x,\chi,q)\exp\Big(D\sqrt{\frac{\log x}{\log\log
      x}}\Big) \ll \log x \underset{\kappa\le k <
  K}{\max}B(T_k,x,\chi,q)\exp\Big(D\sqrt{\frac{\log x}{\log\log
      x}}\Big),\]
so it remains to estimate $B(T_{k},x,\chi,q)$.

To simplify the notation, we write now
 $T_k = T$, $a(T) := (\log\log T)^2$, $b(T) := \frac{\log T}{\log\log
   T}$
and $\mathcal{V}(V,T) := \{n\in\N;\; T\le n <2T,  V_n = V \}$.

We sort the summands corresponding to the values of the $V_{n}$:
\begin{align}
\label{B1}
B(T,x,\chi,q) 
&= \underset{\substack{V\in\N \\ a(T) \le V \le
    b(T)}}{\sum}\underset{\substack{T\le n <2T \\ V_n =
    V}}{\sum}\frac{1}{n}\exp\Big(V\log\Big(\frac{\log x}{\log
      n}\Big) + 2(1 + 2\delta)V\log\log V\Big) \notag\\ 
&\le \frac{1}{T} \underset{\substack{V\in\N \\ a(T) \le V \le
    b(T)}}{\sum}\exp\Big(V\log\Big(\frac{\log x}{\log T}\Big) +
  2(1 + 2\delta)V\log\log V\Big) \operatorname{card}\mathcal{V}(V,T).
\end{align}

Now we split the sum over $V$. For $V\le 2a(T) + 1$ 
we use the trivial estimate
\begin{equation}
\label{card}
  \operatorname{card} \{n\in\N;\; T\le n <2T, ~ V_n = V \} \le T.
\end{equation}
Then we estimate the corresponding sum 
for this part:
\begin{equation}
  \label{B2}
\frac{1}{T} \underset{\substack{V\in\N \\ a(T) \le V \le 2a(T)
    +1}}{\sum}\exp\Big(V\log\Big(\frac{\log x}{\log T}\Big) + 2(1
  + 2\delta)V\log\log V\Big) \operatorname{card}\mathcal{V}(V,T) =
\exp\Big({\rm O}( (\log\log x)^3)\Big).
\end{equation}

Now consider $V\in\N$ with $2a(T)+1 < V \le b(T)$, we split
\[ \mathcal{V}(V,T) = \{n\equiv 0 \text{ mod } 2;\;
  n\in\mathcal{V}(V,T) \} \cup 
  \{n\equiv 1 \text{ mod } 2;\; n\in\mathcal{V}(V,T)
  \}  =:\mathcal{V}_0(V,T) \cup \mathcal{V}_1(V,T).
\] 
Consider a number $n\in\mathcal{V}(V,T)$ for a fixed $V$ 
with $2a(T)+1<V\le b(T)$.
Since $V_n=V$ is the smallest integer such that all $t\in\lbrack n,
n+1\rbrack$ are $V_n$-typical of order $T$, there exists at least one
$t_n\in\lbrack n, n+1 \rbrack$ being $(V_n - 1)$-untypical of order $T$.

So choose for any
$n\in\mathcal{V}(V,T)$ a $t_n\in\lbrack n,n+1 \rbrack$ being
$(V-1)$-untypical. This assignement
gives a bijection between $\mathcal{V}(V,T)$ 
and the set
\[\mathcal{U}(V,T):=\{ t_n ;\;  n\in\mathcal{V}(V,T),
 t_n\in\lbrack n, n+1 \rbrack \text{ and } t_n \text{ is }
 (V-1)\text{-untypical} \} \]
of $(V-1)$-untypical ordinates.
Hence the cardinalities of both sets
are equal, and in $\mathcal{U}(V,T)$ all elements
are $(V-1)$-untypical of order $T$.

Further we define for $h\in\{0,1\}$ the set
\[\mathcal{U}_h(V,T) := \{ t_n\in\mathcal{U}(V,T) ;\;
n\in\mathcal{V}_h(V,T) \}.\]
 
For $t_n\not=t_m$ with $t_n, t_m\in\mathcal{U}_h(V,T)$ 
we have $|t_n-t_m|\ge 1$: If w.l.o.g.\ $n<m$, then 
$t_m-t_n\ge m - (n+1) \ge 1$ since $t_n\in\lbrack n,n+1\rbrack$,
$t_m\in\lbrack m,m+1\rbrack$ and
$n\equiv m$ mod $2$. So the sets
$\mathcal{U}_h(V,T)$ are sets of well distanced 
$(V-1)$-untypical ordinates in the sense
of Proposition \ref{Prop10}.

Since $\operatorname{card}\mathcal{V}(V,T) =
 \operatorname{card}\mathcal{U}(V,T) =
\operatorname{card}\mathcal{U}_0(V,T) 
+ \operatorname{card}\mathcal{U}_1(V,T)$, we can
estimate the cardinality measure of the set
$\mathcal{V}(U,T)$ using Proposition \ref{Prop10},
we obtain
\begin{align}
\label{card1}
\operatorname{card} \mathcal{V}(V,T) &\ll
T\exp\Big(-\frac23(V-1)\log\Big(\frac{V-1}{\log\log T}\Big) +
  \frac43(V-1)\log\log(V-1) 
+ {\rm O}(V) \Big)\notag\\
 &\ll_{\delta} T\exp\Big(-\frac23V\log\Big(\frac{V}{\log\log
       T}\Big)
  +\Big(\frac43+\delta\Big)V\log\log V \Big).  
\end{align}

This leads to the following result:
\begin{align}
 B&(T,x,\chi,q)\le 
 \exp\Big({\rm O}( (\log\log x)^3)\Big) 
\notag\\ 
 &+ \underset{\substack{V\in\N \\ 2a(T) + 1 \le V \le
     b(T)}}{\sum}\frac{1}{T}\exp\Big(V\log\Big(\frac{\log x}{\log
       T}\Big) + 2(1 + 2\delta)V\log\log V\Big)
 \operatorname{card}\mathcal{V}(V,T) 
\text{ by  \eqref{B1} and \eqref{B2} } \notag\\ 
 &\ll_{\delta} \exp\Big({\rm O}(
     (\log\log x)^3)\Big) \notag\\ 
 &+ \underset{\substack{V\in\N \\ 2a(T) + 1 \le V \le b(T)}}{\sum}
 \exp\Big(V\log\Big(\frac{\log x\, (\log\log T)^{2/3}}{\log
       T}\Big) - \frac23 V\log V +
   \Big(\frac{10}{3}+5\delta\Big)V\log\log V \Big) \label{c1}\\ 
  &\ll_{\delta} \exp\Big({\rm O}( (\log\log
      x)^3)\Big) \notag\\ 
  &+ \underset{\substack{V\in\N \\ 2a(T) + 1 \le V \le b(T)}}{\sum}
  \exp\Big(V\log\Big(\frac{\log x \, \log\log T}{\log T}\Big) -
    \frac23 V\log V + \Big(\frac{10}{3}+5\delta\Big)V\log\log V
  \Big), \label{c2} 
\end{align}
where in \eqref{c1} the implicit constant in the estimate
depends on $\delta$ since we used equation \eqref{card1}.

In order to majorize the last sum \eqref{c2}, we use Proposition
\ref{Prop20} with the following parameters:
\[ \text{Let } A := \log\Big(\frac{\log x \,\log\log T}{\log T}\Big)
 \text{ and } C := \frac{10}{3} + 5\delta. \]
(Then $A\ge 4C^4 + 1$ and $V>e^{3C/2}$ hold if $x$ is large
enough.) 

We obtain
\begin{align}
&\underset{\substack{V\in\N \\ 2a(T) + 1 \le V \le b(T)}}{\sum}
\exp\Big(V\log\Big(\frac{\log x \,\log\log T}{\log T}\Big) -
  \frac23 V\log
  V + \Big(\frac{10}{3}+5\delta\Big)V\log\log V \Big)\notag\\ 
&\le ~ \frac{\log T}{\log\log T}\exp\Big(\Big(\log x
    \frac{\log\log T}{\log T}\Big)^{3/2}\Big(\frac32 \log\Big(\log
      x\frac{\log\log T}{\log T}\Big)\Big)^{5 +
    15\delta/2} \Big). \label{c3} 
\end{align}
Since
\begin{equation*}
\frac{\log\log T}{\log T} = \frac{\log\log T_k}{\log T_k} \le
\frac{\log\log T_{\kappa}}{\log T_{\kappa}}\ll \frac{\log\log x}{(\log
  x)^{3/5}(\log\log x)^c} \le (\log x)^{-3/5},
\end{equation*}
we have
\[\Big(\log x \frac{\log\log T}{\log T}\Big)^{3/2}
\le \Big((\log x)^{2/5}
  (\log\log x)^{1-c}\Big)^{3/2} = (\log x)^{3/5} (\log\log
x)^{3/2-3c/2},\]
and as $c\ge 1$, we obtain further
\begin{equation*}
\Big(\frac32 \log\Big(\log x\frac{\log\log T}{\log
      T}\Big)\Big)^{5 + 15\delta/2}
 \le(\log\log x)^{5 + 15\delta/2}.
\end{equation*}

Using these estimates, we continue
the estimation of \eqref{c3} with
\begin{align*}
&\le \exp\Big( \log\log x + (\log x)^{3/5} 
  (\log\log x)^{3/2-3c/2+5+15\delta/2} \Big)\\
&=\exp\Big((\log x)^{3/5}(\log\log x)^{13/2 - 3c/2
    + 15\delta/2} + \log\log x \Big) \\ 
&\ll_{\delta} \exp\Big((\log x)^{3/5}(\log\log x)^{13/2 -
    3c/2 + 8\delta}\Big).
\end{align*}

Now we resume everything including the
term $\exp\Big(D\sqrt{\frac{\log x}{\log\log x}}\Big)$ again, 
we obtain
\[B(x,\chi,q) \ll_{\delta} \exp\Big((\log x)^{3/5}(\log\log
  x)^{13/2 - 3c/2 + 8\delta}
\Big)\exp\Big((D+1)\sqrt{\frac{\log x}{\log\log x}}\Big), \]
and using the estimate
\begin{align*}
(\log x)^{3/5}&(\log\log x)^{13/2 - 3c/2 + 8\delta}
+ (D+1) \sqrt{\frac{\log x}{\log\log x}}\\
 &\ll (\log
x)^{3/5}(\log\log x)^{13/2 - 3c/2 + 8\delta}
\Big(1 + \log(x)^{-1/10}(\log\log x)^{3c/2}\Big) \\ 
&\ll (\log x)^{\frac{3}{5}}(\log\log x)^{13/2 - 3c/2 + 8\delta},
\end{align*}
we obtain finally
\[B(x,\chi,q) \ll_{\delta} \exp\Big((\log x)^{3/5}(\log\log
  x)^{13/2 - 3c/2 + 8\delta} \Big).\]
\qed

\Leer
Now we still have to consider the principal character mod $q$,
for this we use the result of the zeta-function.

\begin{lemma}
\label{l3} 
Let $q\in\N$, $x\ge q>1$, then we have for the principal character
$\chi_0$ mod $q$ the estimate
\[ \text{A}(x,\chi_0,q) \ll_{\delta} \sqrt{x} \exp\Big((\log
    x)^{1/2} (\log\log x)^{5/2 + 4\delta}\Big). \] 
\end{lemma}

\textbf{Proof:}
Due to the formula
\[\Ls(s,\chi_0) = \zeta(s)\underset{p|q}{\prod}\Big(1-\frac{1}{p^s}\Big), \]
we use the estimate for the zeta-integral.
So we estimate the product
$\Big|\prod_{p|q}(1-p^{-s})^{-1}\Big|$ for $\sigma\ge\frac12$.

For this, consider the logarithm of the product and include the series
expansion of the logarithm:
\begin{align*}
\Big|\underset{p|q}{\sum}-\log\Big(1-\frac{1}{p^s}\Big)\Big| &=
  \Big|\underset{p|q}{\sum}-\underset{k\in\N}{\sum}
  (-1)^{k+1}\frac{(-p^{-s})^k}{k}\Big| =
  \Big|\underset{p|q}{\sum}\underset{k\in\N}{\sum}
  (-1)^{2k+2}\frac{1}{kp^{ks}}\Big|\\ 
&\le \underset{p|q}{\sum}\underset{k\in\N}{\sum}\frac{1}{kp^{k/2}}
    = \underset{p|q}{\sum}\frac{1}{p^{1/2}} + \frac12
\underset{p|q}{\sum}\frac1p +
\underset{p|q}{\sum}\underset{k>2}{\sum}\frac{1}{kp^{k/2}} \\ 
&\le \underset{i=1}{\overset{2\log q}{\sum}}\frac{1}{p_i^{1/2}} +
  \frac{1}{2}\underset{p\le q}{\sum}\frac{1}{p} + {\rm O}(1) \\
&\ll \sqrt{\frac{\log q}{\log\log q}} + \log\log q +  {\rm O}(1).
\end{align*}	

We conclude
\[|\Ls(s,\chi_0)|^{-1}\ll|\zeta(s)|^{-1}\exp\Big(D\sqrt{\frac{\log
      q}{\log\log q}}\Big) \] 
for an absolute constant $D>0$.

Since $\sqrt{\frac{\log q}{\log\log q}}$ is
monotonic increasing in $q$, we have for $x \ge q$ 
\[\Ls(s,\chi_0)^{-1}\ll\zeta(s)^{-1}\exp\Big(D\sqrt{\frac{\log
      x}{\log\log x}}\Big).\]

Now the additional term $\sqrt{\frac{\log
x}{\log\log x}}$ does not disturb the magnitude of the exponent
in the final result, since we have
\begin{align*}
\Big|\underset{S(x,\chi,q)}{\int}\Ls(z,\chi_0)^{-1}\frac{x^z}{z}dz \Big|
&\ll \underset{S(x,\chi,q)}{\int}\Big|\zeta(z)^{-1}\frac{x^z}{z}\Big|dz
   \exp\Big(D\sqrt{\frac{\log x}{\log\log x}}\Big) \\
&\ll_{\delta}\sqrt{x}\exp\Big((\log x)^{1/2}(\log\log
  x)^{5/2+ 4\delta} 
   + D\sqrt{\frac{\log x}{\log\log x}}\Big)\\ 
&\ll\sqrt{x}\exp\Big((\log x)^{1/2}(\log\log x)^{5/2+ 4\delta}\Big), 
\end{align*}
where we have set $c=\frac52 + 3\delta$ in the estimate 
at the end of the paper of \cite{RoBa}.
\qed

\Leer
\textbf{Proof of Theorem \ref{t}:}

Let $q > 2$, since for $q=2$ there is only the principal character and
we can use then the sharper result from Lemma \ref{l3}.

We use equation \eqref{vorueb}, Lemma \ref{l1} and Lemma \ref{l2} and
set $c=\frac{11}{5} + \frac{16}{5}\delta$, together with
Lemma \ref{l3} we obtain
\begin{align*}
\Big|M(x,a,q)\Big|&\le
\frac{1}{\varphi(q)}\underset{\chi (q)}{\sum}\Big|\underset{n\le
    x}{\sum}\chi(n)\mu(n)\Big| =
\frac{1}{\varphi(q)}\underset{\chi (q)}{\sum}\Big|A(x,\chi,q)\Big|
+ {\rm O}(\log x) \\ 
&=\frac{1}{\varphi(q)}|A(x,\chi_0,q)| +
\frac{1}{\varphi(q)}\underset{\substack{\chi (q)\\
    \chi\not=\chi_0}}{\sum}\Big|A(x,\chi,q)\Big| + {\rm O}(\log x) \\ 
&\ll_{\delta} \frac{1}{\varphi(q)} \sqrt{x}\exp\Big((\log x)^{1/2} 
  (\log\log x)^{5/2 +  4\delta}\Big) \\
&+ \frac{\varphi(q) -1}{\varphi(q)}\sqrt{x}\exp\Big((\log x)^{3/5} 
   (\log\log x)^{16/5 + 16\delta/5}\Big) \\ 
&\ll\sqrt{x}\exp\Big((\log x)^{3/5} (\log\log x)^{16/5 + 16\delta/5}\Big).
\end{align*}

Since $\delta\in(0,1\rbrack$ can be choosen arbitrary, we get the
assertion with the choice $\delta = \frac{5}{16}\varepsilon$.
\qed



\vspace{\baselineskip}

Karin Halupczok\\
Westf\"{a}lische Wilhelms-Universit\"{a}t M\"{u}nster\\
Mathematisches Institut\\
Einsteinstra\ss{}e 62 \\
D-48149 M\"{u}nster \\
\texttt{karin.halupczok@uni-muenster.de}

\vspace{\baselineskip}

Benjamin Suger\\
Albert-Ludwigs-Universit\"{a}t Freiburg\\
Department of Computer Science\\
Autonomous Intelligent Systems Group\\
Georges-Koehler-Allee 079\\
D-79110 Freiburg, Germany\\
\texttt{suger@informatik.uni-freiburg.de}

\end{document}